\documentclass[12pt]{amsart}

\copyrightinfo{2000}{American Mathematical Society}

\newcommand{\cD}{{\mathcal D}}

\newcommand{\cP}{{\mathcal P}}
\newcommand{\cS}{{\mathcal S}}

\newcommand{\bN}{{\mathbb N}}
\newcommand{\bZ}{{\mathbb Z}}
\newcommand{\bQ}{{\mathbb Q}}
\newcommand{\bR}{{\mathbb R}}
\newcommand{\bC}{{\mathbb C}}

\numberwithin{equation}{section}

\newtheorem{Theorem}{Theorem}[section]
\newtheorem{Lemma}{Lemma}[section]

\newtheorem{Definition}{Definition}[section]
\newtheorem{Remark}{Remark}[section]

\newtheorem{Proposition}{Proposition}[section]

\author{V.~M.~Shelkovich}
\address{Department of Mathematics, St.-Petersburg State Architecture
and Civil Engineering University, \ 2 Krasnoarmeiskaya 4, 190005,
St. Petersburg, \ Russia.}
\email{shelkv@vs1567.spb.edu}

\title[Tauberian theorems for distributions in the Lizorkin spaces]
{Tauberian theorems for distributions in the Lizorkin spaces}

\thanks{This paper was supported in part by DFG Project 436 RUS 113/809/0-1.}

\subjclass[2000]{Primary 40E05, 46F10; Secondary 26A33, 46F12}

\keywords{Tauberian type theorems, the Lizorkin spaces, distributions,
fractional operator}

\date{ }

\begin{document}

\begin{abstract}
In this paper some multidimensional Tauberian theorems for the Lizorkin
distributions (without restriction on the support) are proved.
Tauberian theorems of this type are connected with the Riesz
fractional operators.
\end{abstract}

\maketitle

\section{Introduction}
\label{s1}

\subsection{Tauberian type theorems.}\label{s1.1}
As is well known, in mathematical physics there are so-called
{\it Tauberian theorems} which have many applications. Tauberian
theorems are usually assumed to connect asymptotical behavior of a
function (distribution) at zero with asymptotical
behavior of its Fourier (Laplace or other integral transforms) at infinity.
The inverse theorems are usually called ``Abelian''~\cite{D-Zav1},
~\cite{D-Zav2},~\cite{Kor},~\cite{Vl-D-Zav} (see also the
references cited therein).

Multidimensional Tauberian theorems for distributions (as a rule,
from the space $\cS'(\bR^n)$) have been treated by V.~S.~Vladimirov,
Yu.~N.~Drozzinov, B.~I.~Zavyalov in the fundamental book~\cite{Vl-D-Zav}.
It remains to note that  Tauberian theorems have many applications and
are intensively used in mathematical physics (see~\cite{Kor},~\cite{Vl-D-Zav}).

Some types of Tauberian theorems are connected with the
fractional operator. In~\cite{Vl-D-Zav}, as a rule, theorems of this type
were proved for distributions with supports in
the cone in $\bR^n$, $n\ne 1$ or in semiaxis for $n=1$. This is related
to the fact that such distributions constitute a {\it convolution algebra\/}.
In this case a kernel of the fractional operator is the distribution
with a support in a cone in $\bR^n$, $n\ne 1$ or in semiaxis for
$n=1$~\cite[\S2.8.]{Vl-D-Zav}. Thus in this case, in general, the
convolution of a distribution and a kernel of the fractional operator
is not well defined in the sense of the space $\cS'(\bR^n)$). Moreover,
in general, the Schwartian test function space ${\cS}(\bR^n)$
{\it is not invariant\/} under the fractional operators.
In view of this fact Tauberian type theorems
for distributions {\it without restriction on the support\/} have
not been considered in~\cite{Vl-D-Zav}.

The solution of the above problem was suggested by P.~I.~Lizorkin in the
excellent papers~\cite{Liz1}--~\cite{Liz3}(see also~\cite{Sam1},~\cite{Sam2}).
Namely, in~\cite{Liz1}--~\cite{Liz3} a new type spaces
{\it invariant\/} under fractional operators were introduced (see
Lemmas~\ref{lem2},~\ref{lem2.1}).
The Lizorkin spaces are ``natural'' definitional domains of the fractional
operators.
Note that fractional operators have many applications and are intensively
used in mathematical physics~\cite{Eid-Koch}, \cite{Sam3},
~\cite{Sam-Kil-Mar}. These two last fundamental books have the exhaustive
references.

Thus, if we  want to prove Tauberian type theorems for distributions
{\it without restriction on the supports\/}, we must consider distributions
from the Lizorkin spaces.

In this paper some multidimensional Tauberian theorems for
the Lizorkin distributions ({\it without restriction on the support\/})
are proved. Tauberian theorems of this type are connected with
the Riesz fractional operator.

This paper was inspired by our paper on $p$-adic Tauberian
theorems~\cite{Kh-S1}. The point is, that in $p$-adic analysis
a kernel $f_{\alpha}(z)=\frac{|z|_p^{\alpha-1}}{\Gamma_p(\alpha)}$
of the $p$-adic Vladimirov fractional operator
$D^{\alpha}=f_{-\alpha}*$ is the $p$-adic a distribution
{\it without restriction on the support\/}, where $*$ is a
convolution, $|z|_p^{\alpha-1}$ is a $p$-adic homogeneous
distribution, $\Gamma_p(\alpha)$ is the $p$-adic $\Gamma$-function
(see~\cite{Vl-V-Z}).

\subsection{Contents of the paper.}\label{s1.2}
In Sec.~\ref{s2} we recall some facts from the theory of distributions.
In particular, in Subsec.~\ref{s2.1} some properties of the
{\it Lizorkin spaces\/} of test functions and
distributions~\cite{Liz1}--~\cite{Liz3},~\cite{Sam3},
~\cite{Sam-Kil-Mar} are given.
In Subsec.~\ref{s2.2} we recall properties of the
{\it Riesz potential}~\cite{Sam3},~\cite{Sam-Kil-Mar}.
In Subsec.~\ref{s2.3},~\ref{s2.3-1} the {\it Riesz fractional operator},
which was studied in~\cite{Sam3},~\cite{Sam-Kil-Mar} is
introduced. In Subsec.~\ref{s2.4} by Definition~\ref{de4.1}
we give the notion of regular variation introduced by J.~Karamata.
In this subsection Definitions~\ref{de4},~\ref{de5} of the
{\it quasi-asymptotics\/} at infinity and at zero for
distributions~\cite{D-Zav1},~\cite{Vl-D-Zav} are introduced.

In Sec.~\ref{s3}, some multidimensional Tauberian type theorems
(Theorems~\ref{th5}--~\ref{th9}) for distributions are proved.
Theorems~\ref{th5},~\ref{th6} are
related to the Fourier transform and hold for distributions from
${\cS}'(\bR^n)$. Theorems~\ref{th7}--~\ref{th9} are related to the
fractional operators and hold for distributions from
the Lizorkin spaces $\Phi_{\times}'(\bR^n)$ and $\Phi'(\bR^n)$.

\section{Some results from the theory of distributions}
\label{s2}

\subsection{The Lizorkin spaces.}\label{s2.1}
We denote by $\bN$, $\bZ$, $\bR$, $\bC$ the sets of positive integers,
integers, real, complex numbers respectively, and set ${\bN}_0={0}\cup{\bN}$.
If $x=(x_1,\dots,x_n)$ then $|x|=\sqrt{x_1^2+\cdots+x_n^2}$
and $x^{j}\stackrel{def}{=}x_1^{j_1}\cdots x_n^{j_n}$.
For $j=(j_1,\dots,j_n)\in {\bN}_0^n$ we assume $j!=j_1!\cdots j_n!$,
\ $|j|=j_1+\cdots+j_n$.
We shall denote partial derivatives of the order $|j|$ by
$\partial_x^{j}=\frac{\partial^{|j|}}
{\partial{x_1}^{j_1}\cdots\partial{x_n}^{j_n}}$.

Denote by ${\cD}(\bR^n)$ and ${\cS}(\bR^n)$ the linear spaces of
infinitely differentiable functions with a compact support
and the Schwartian test function space.
Denote by ${\cD}'(\bR^n)$ and ${\cS}'(\bR^n)$ the space of all linear
and continuous functionals on ${\cD}(\bR^n)$ and ${\cS}(\bR^n)$,
respectively (see~\cite{Brem},~\cite{Vl-D-Zav}).

\begin{Definition}
\label{de2} \rm
(~\cite[Ch.III,\S3.1.]{Gel-S}) A distribution $f \in {\cS}'(\bR^n)$
is called {\it homogeneous\/} of degree $\alpha$ if
$$
\Bigl\langle f,\varphi\Big(\frac{x_1}{t},\dots,\frac{x_n}{t}\Big)\Bigr\rangle
=t^{\alpha+n}\langle f,\varphi \rangle,
\quad \forall \, \varphi \in {\cS(\bR^n)}, \quad t>0,
$$
($t\in \bR$), i.e.,
$$
f(tx_1,\dots,tx_n)=t^{\alpha}f(x_1,\dots,x_n), \quad \forall \, t>0.
$$
\end{Definition}

By reformulating our definition~\cite{Al-Kh-S1},~\cite{Al-Kh-S2}
for the case of $\bR^n$ (instead of the field of $p$-adic numbers $\bQ_p$),
we introduce the following definition.

\begin{Definition}
\label{de2.1} \rm
A distribution $f_m \in {\cS}'(\bR^n)$ is said to be
{\it associated homogeneous {\rm(}in the wide sense{\rm)}\/}
of degree~$\alpha$ and order~$m$, \ $m=0,1,2,\dots$, if
$$
\Bigl\langle f_m,\varphi\Big(\frac{x_1}{t},\dots,\frac{x_n}{t}\Big)\Bigr\rangle
=t^{\alpha+n}\langle f_m,\varphi \rangle,
+\sum_{j=1}^{m}t^{\alpha+n}\log^jt\langle f_{m-j},\varphi \rangle,
$$
for all $\varphi \in {\cS(\bR^n)}$ and $t>0$ ($t\in \bR$), where
$f_{m-j}$ is an {\it associated homogeneous {\rm(}in the broad
sense{\rm)}\/} distribution of degree~$\alpha$ and order $m-j$, \
$j=1,2,\dots,m$, i.e.,
$$
f_m(tx_1,\dots,tx_n)=t^{\alpha}f_m(x_1,\dots,x_n)
+\sum_{j=1}^{m}t^{\alpha}\log^jtf_{m-j}(x_1,\dots,x_n),
\quad \forall \, t>0.
$$
Here for $m=0$ the sum is empty.
\end{Definition}

Associated homogeneous distributions (in the wide sense) of order
$m=1$ coincide with associated homogeneous distributions of order $m=1$.
Associated homogeneous distributions of order $m=0$ coincide with
homogeneous distributions.

\begin{Remark}
\label{rem1} \rm
We recall that the notion of the {\it associated homogeneous distribution\/}
from ${\cD}'(\bR)$ was introduced in~\cite[Ch.I,\S 4.1.]{Gel-S}
by the following definition:
for any~$m$, the distribution $f_m\in {\cD}'(\bR)$ is called
an {\it associated homogeneous distribution of order~$m$ and
degree\/}~$\alpha$ if for any $t >0$ and any
$\varphi \in {\cD}(\bR)$ we have
\begin{equation}
\label{3}
\Bigl\langle f_m,\varphi\Big(\frac{x}{t}\Big) \Bigr\rangle
=t^{\alpha+1} \langle f_m,\varphi \rangle
+ t^{\alpha+1}\log{t}\langle f_{m-1},\varphi \rangle,
\end{equation}
where $f_{m-1}$ is an {\it associated homogeneous distribution\/}
of order $m-1$ and of degree~$\alpha$, \ $m=1,2,3,\dots$. In
the paper~\cite[Ch.X,8.]{Vil}, giving a brief outline of the
book~\cite{Gel-S}, the definition of an {\it associated homogeneous
distribution\/} was introduced as an analog of relation (\ref{3}),
where in the right-hand side of (\ref{3}) $\log{t}$ is replaced by
$\log^m{t}$. Definition~(\ref{3}) is introduced by analogy with
the definition of an associated eigenvector.

In the book~\cite[Ch.I,\S 4.2.]{Gel-S} and in the
paper~\cite[Ch.X, 8.]{Vil} it is stated that the distributions
$x_{\pm}^{\alpha}\log^m x_{\pm}$, $\alpha \ne -1,-2,\dots$
and $P\big(x_{\pm}^{-n}\log^{m-1}x_{\pm}\big)$ are
{\it associated homogeneous distributions\/} of order $m$ and
degree~$\alpha$ and $-n$, respectively, $m=1,2,3,\dots$.
If $m=2,3,\dots$, it is easily verified that these distributions
are not associated homogeneous in the sense of Definition~(\ref{3})
and the {\it modified definition\/} from~\cite[Ch.X, 8.]{Vil}.
We illustrate this fact by the following simple example:
$$
\log^2(tx_{\pm})=\log^2 x_{\pm}+2\log t \log x_{\pm}+\log^2t,
\quad t>0.
$$
One can prove that associated homogeneous (in {\it the strict sense\/} )
distributions can only have order~$1$, while for $m\ge 2$
Definition~(\ref{3}) describes an empty class.

Thus an associated homogeneous distribution ({\it in the wide sense\/})
$f_m$ of order $m$, $m\ge 2$ is reproduced by the similitude operator
$U_{a}f(x)=f(ax)$ up to a linear combination of associated homogeneous
distributions ({\it in the wide sense\/}) of orders
$m-1,m-2,\dots,0$, and therefore, strictly speaking, it is not
an {\it associated homogeneous distribution\/}.
Following the book~\cite[Ch.I,\S 4.1.]{Gel-S}, even for $m\ge 2$, we will
call these distributions a.h.d., omitting the words ``in the wide sense''.
One can prove that ($\bC$) distributions $x_{\pm}^{\alpha}\log^m x_{\pm}$,
$\alpha \ne -1,-2,\dots$, and $P\big(x_{\pm}^{-n}\log^{m-1}x_{\pm}\big)$
are associated homogeneous distributions ({\it in the wide wide sense\/}) in
terms of Definition~\ref{de2.1} for the case $n=1$.

The above mentioned problem is not considered in the present paper.
\end{Remark}

The Fourier transform of $\varphi\in {\cS}(\bR^n)$ is defined by the
formula
$$
F[\varphi](\xi)=\int_{\bR^n}e^{i\xi\cdot x}\varphi(x)\,d^nx,
\quad \xi \in \bR^n,
$$
where $\xi\cdot x$ is the scalar product of vectors. It is well known
that the Fourier transform is a linear isomorphism ${\cS}(\bR^n)$ into
${\cS}(\bR^n)$.
We define the Fourier transform $F[f]$ of a distribution
$f\in{\cS}'(\bR^n)$ by the relation
\begin{equation}
\label{5.2}
\langle F[f],\varphi\rangle=\langle f,F[\varphi]\rangle,
\end{equation}
for all $\varphi\in {\cS}(\bR^n)$.

For distributions $f,g\in{\cS}'(\bR^n)$ the convolution $f*g$ is
defined as
\begin{equation}
\label{6}
\langle f*g,\varphi\rangle=\langle f(x)\times g(y),\varphi(x+y)\rangle,
\end{equation}
for all $\varphi\in {\cS}(\bR^n)$, where $f(x)\times g(y)$ is
the direct product of distributions.
If for distributions $f,g\in{\cS}'(\bR^n)$ a convolution $f*g$
exists then
\begin{equation}
\label{7}
F[f*g]=F[f]F[g].
\end{equation}

Recall the well known facts from~\cite{Liz1}--~\cite{Liz3},~\cite[2.]{Sam3},
~\cite[\S 25.1.]{Sam-Kil-Mar}. Consider the following subspace of
the space $\cS(\bR^n)$
$$
\Psi=\Psi(\bR^n)=\{\psi(\xi)\in \cS(\bR^n):
(\partial_{\xi}^{j}\psi)(0)=0, \, |j|=1,2,\dots\}.
$$
Obviously, $\Psi \ne \emptyset$. The space of functions
$$
\Phi=\Phi(\bR^n)=\{\phi: \phi=F[\psi], \, \psi\in \Psi(\bR^n)\}.
$$
is called the {\it Lizorkin space of test functions\/}.
The Lizorkin space can be equipped with the topology of the space $\cS(\bR^n)$
which makes $\Phi$ a complete space~\cite[2.2.]{Sam3},
~\cite[\S 25.1.]{Sam-Kil-Mar}.

Since the Fourier transform is a linear isomorphism ${\cS}(\bR^n)$ into
${\cS}(\bR^n)$, this space admits the following characterization:
$\phi\in \Phi$ if and only if $\phi\in \cS(\bR^n)$ is orthogonal to
polynomials, i.e.,
\begin{equation}
\label{8}
\int_{\bR^n}x^{j}\phi(x)\,d^nx=0, \quad |j|=0,1,2,\dots.
\end{equation}

Let $\cP \subset \cS'$ be the subspace of all polynomials. As is well
known, the set $F[\cP]\subset \cS'$ of the Fourier transform of polynomials
consists of finite linear combinations of the Dirac $\delta$-function
supported at the origin and its derivatives.
Here $\cP=\Phi^{\perp}$ and $F[\cP]=\Psi^{\perp}$, where $\Phi^{\perp}$
and $\Psi^{\perp}$ are the subspaces of functionals in $\cS'$ which
are orthogonal to $\Phi$ and $\Psi$, respectively~\cite[2.]{Sam3},
~\cite[\S 8.2.]{Sam-Kil-Mar}.

\begin{Proposition}
\label{pr1}
{\rm (~\cite[Proposition~2.5.]{Sam3})} The spaces of linear and
continuous functionals $\Phi'$ and $\Psi'$ can be identified with
the quotient spaces
$$
\Phi'=\cS'/\cP, \qquad
\Psi'=\cS'/F[\cP]
$$
modulo the subspaces $\cP$ and $F[\cP]$, respectively.
\end{Proposition}

The space $\Phi'$ is called the {\it Lizorkin space of distributions\/}.

Analogously to (\ref{5.2}), we define the Fourier transform of
distributions $f\in \Phi'(\bR^n)$ and $g\in \Psi'(\bR^n)$ by the
relations~\cite[(25.18),(25.18')]{Sam-Kil-Mar}:
\begin{equation}
\label{8.1}
\begin{array}{rcl}
\displaystyle
\langle F[f],\psi\rangle=\langle f,F[\psi]\rangle,
&& \forall \, \psi\in \Psi(\bR^n), \medskip \\
\displaystyle
\langle F[g],\phi\rangle=\langle g,F[\phi]\rangle,
&& \forall \, \phi\in \Phi(\bR^n). \\
\end{array}
\end{equation}
By definition, $F[\Phi(\bR^n)]=\Psi(\bR^n)$ and
$F[\Psi(\bR^n)]=\Phi(\bR^n)$, i.e., these definitions are correct.

Now we introduce another type of the Lizorkin space.
Let
$$
\Psi_{\times}=\Psi_{\times}(\bR^n)=\bigl\{\psi(\xi)\in \cS(\bR^n):
\qquad\qquad\qquad\qquad\qquad\qquad\qquad\qquad
$$
$$
\qquad
(\partial_{\xi}^{j}\psi)(\xi_1,\dots,\xi_{k-1},0,\xi_{k+1},\dots,\xi_{n})=0,
\, |j|=1,2,\dots, \, k=1,2,\dots,n\bigr\}.
$$
The space of functions
$$
\Phi_{\times}=\Phi_{\times}(\bR^n)=\{\phi: \phi=F[\psi],
\, \psi\in \Psi(\bR^n)\}.
$$
is called the {\it Lizorkin space of test functions\/}.
This space can be equipped with the topology of the space
$\cS(\bR^n)$ which makes $\Phi_{\times}$ a complete space.

It is clear that $\phi\in \Phi_{\times}$ if and only if
\begin{equation}
\label{8*}
\int_{-\infty}^{\infty}x^{m}_{k}
\phi(x_1,\dots,x_{k-1},x_{j},x_{k+1},\dots,x_{n})\,dx_k=0,
\end{equation}
$m=0,1,2,\dots$, \ $k=1,2,\dots,n$.

Analogously to Proposition~\ref{pr1},
$$
\Phi'_{\times}=\cS'/\Phi^{\perp}_{\times}, \qquad
\Psi'_{\times}=\cS'/\Psi^{\perp}_{\times},
$$
where $\Phi^{\perp}_{\times}$ and $\Psi^{\perp}_{\times}$ are
subspaces of functionals in $\cS'$ which are orthogonal to $\Phi_{\times}$
and $\Psi_{\times}$, respectively.

We define the Fourier transform of
distributions $f\in \Phi'_{\times}(\bR^n)$ and $g\in \Psi'_{\times}(\bR^n)$
similarly to Definition (\ref{8.1}).

\subsection{The Riesz potentials.}\label{s2.2}
Let us introduce the distribution $|x|^{\alpha}\in {\cS}'(\bR^n)$
(see~\cite[Lemma~2.9.]{Sam3},~\cite[(25.19)]{Sam-Kil-Mar},
~\cite[Ch.I,\S3.9.]{Gel-S}).
If $Re\,\alpha>-n$ then the function $|x|^{\alpha}$ is locally
integrable and generates a regular functional
\begin{equation}
\label{9}
\langle |x|^{\alpha},\varphi(x)\rangle
=\int_{\bR^n}|x|^{\alpha}\varphi(x)\,d^nx,
\quad \forall \, \varphi\in {\cS}(\bR^n).
\end{equation}
If $Re\,\alpha \le -n$, we define this distribution by means of
analytic continuation:
$$
\langle |x|^{\alpha},\varphi \rangle
=\int_{|x|<1}|x|^{\alpha}\Big(\varphi(x)
-\sum_{|j|=0}^{m}\frac{x^{j}}{j!}(\partial_x^{j}\varphi)(0)\Big)\,d^nx
\qquad\qquad\qquad\qquad\qquad
$$
\begin{equation}
\label{10}
+\int_{|x|>1}|x|^{\alpha}\varphi(x)\,d^nx
+\sum_{k=0}^{[m/2]}\frac{\pi^{\frac{n}{2}}(\Delta^{k}\varphi)(0)}
{2^{2k-1}k!\Gamma(\frac{n}{2}+k)(\alpha+n+2k)},
\end{equation}
for all $\varphi\in {\cS}(\bR^n)$, where $\alpha+n\ne 0,-2,-4,\dots$,
and $m >-Re\,\alpha -n-1$, \ $\Delta$ is the Laplacian,
$[a]$ is the integral part of a number $a$. Relation
(\ref{10}) gives the explicit formula of analytic continuation of the
distribution $|x|^{\alpha}$ from the half-plain $Re\,\alpha>-n$ to the
domain $-m-n-1< Re\,\alpha \le -n$.

Formula (\ref{10}) is proved by using the relation~\cite[\S25.1]{Sam-Kil-Mar}
$$
\sum_{|j|=0}^{m}\frac{(\partial_x^{j}\varphi)(0)}{j!}
\int_{|x|<1}x^{j}|x|^{\alpha}\,d^nx
=\sum_{k=0}^{[m/2]}\frac{\pi^{\frac{n}{2}}(\Delta^{k}\varphi)(0)}
{2^{2k-1}k!\Gamma(\frac{n}{2}+k)(\alpha+n+2k)},
$$
and formula $\Delta^{m}=\sum_{|j|=m}\frac{m!}{j!}\partial_x^{2j}$,
\ $m=0,1,2,\dots$.

It is clear that the distribution $|x|^{\alpha}$, \ $\alpha \ne -n-2s$,
$s \in {\bN}_0$ is a homogeneous distribution of degree~$\alpha$
(see Definition~\ref{de2}).

In the case $\alpha \ne -n-2s$, $s\in {\bN}_0$, excluded in (\ref{10}),
according to~\cite[(2.29)]{Sam3},~\cite[(25.23),(25.24)]{Sam-Kil-Mar},
we define $\langle |x|^{\alpha},\varphi(x)\rangle$ as
$$
\Bigl\langle P\Big(\frac{1}{|x|^{n+2s}}\Big),\varphi \Bigr\rangle
\qquad\qquad\qquad\qquad\qquad\qquad\qquad\qquad\qquad\qquad\qquad
$$
\begin{equation}
\label{12}
=\lim_{\alpha\to -n-2s}\bigg(\langle |x|^{\alpha},\varphi \rangle
-\frac{\pi^{\frac{n}{2}}(\Delta^{s}\varphi)(0)}
{2^{2s-1}s!\Gamma(\frac{n}{2}+s)(\alpha+n+2s)}\bigg),
\end{equation}
for all $\varphi\in {\cS}(\bR^n)$. Here the distribution
$P(|x|^{-n-2s})$ is called the principal value of the
function~$|x|^{-n-2s}$. In view of Definition~\ref{de2.1},
this distribution is an {\it associated homogeneous\/} distribution
of degree $-n-2s$ and order $1$.

Thus the distribution $|x|^{\alpha}$ is defined for any $\alpha \in \bC$.

Let us introduce the distribution from ${\cS}'(\bR^n)$
\begin{equation}
\label{15}
\kappa_{\alpha}(x)=\frac{|x|^{\alpha-n}}{\gamma_n(\alpha)},
\quad \alpha \ne -2s, \, \alpha \ne n+2s, \quad s=0,1,2,\dots
\end{equation}
called the {\it Riesz kernel\/}, where $|x|^{\alpha}$ is a
homogeneous distribution of degree~$\alpha$ defined by (\ref{10}),
\begin{equation}
\label{15.1}
\gamma_n(\alpha)=\frac{2^{\alpha}\pi^{\frac{n}{2}}\Gamma(\frac{\alpha}{2})}
{\Gamma(\frac{n-\alpha}{2})}.
\end{equation}
The Riesz kernel is an entire function of the complex variable $\alpha$.

In view of~\cite[Ch.I,\S3.9.,(8)]{Gel-S},~\cite[4.,(69),(71)]{Kan},
$$
\frac{|x|^{\alpha-n}\Gamma(\frac{n}{2})}{\pi^{\frac{n}{2}}
\Gamma(\frac{\alpha}{2})}
\biggr|_{\alpha=-2s}
\stackrel{def}{=}\lim_{\alpha \to -2s}
\frac{|x|^{\alpha-n}\Gamma(\frac{n}{2})}{\pi^{\frac{n}{2}}
\Gamma(\frac{\alpha}{2})}
=\frac{{\rm Res}_{\alpha=-2s}|x|^{\alpha-n}}
{{\rm Res}_{\alpha=-2s}\langle |x|^{\alpha-n},e^{-|x|^2}\rangle }
\qquad
$$
\begin{equation}
\label{16}
=\left\{
\begin{array}{rcl}
\displaystyle
\delta(x), \quad  s&=&0, \smallskip \\
\displaystyle
\frac{(-1)^{s}}{2^{s}n(n+2)\cdots (n+2s-2)}\Delta^{s}\delta(x),
\quad s&=&1,2,\dots, \\
\end{array}
\right.
\end{equation}
where ${\rm Res}$ stands for the residue, and the limit is understood
in the weak sense.
Using formulas (\ref{16}), (\ref{15}), (\ref{15.1}), and
$\Gamma(\frac{n}{2}+s)=2^{-s}n(n+2)\cdots (n+2s-2)\Gamma(\frac{n}{2})$,
we define $\kappa_{-2s}(x)$ as a distribution from ${\cS}'(\bR^n)$:
\begin{equation}
\label{17}
\kappa_{-2s}(x)\stackrel{def}{=}\lim_{\alpha \to -2s}\kappa_{\alpha}(x)
=(-\Delta)^{s}\delta(x), \quad s=0,1,2,\dots,
\end{equation}
where the limit is understood in the weak sense.

Next, using formulas (\ref{9}), (\ref{15}), (\ref{15.1}), and formula
$$
\Gamma\Big(-\frac{\beta}{2}-s\Big)
=\frac{(-1)^{s+1}2^{s+1}\Gamma(-\frac{\beta}{2}+1)}
{\beta(\beta+2)\cdots(\beta+2s-2)(\beta+2s)},
$$
and taking into account (\ref{8}), we define $\kappa_{n+2s}(x)$,
$s=0,1,2,\dots$ as distribution from the {\it Lizorkin space of
distributions\/} $\Phi'(\bR^n)$:
$$
\langle \kappa_{n+2s}(x),\phi \rangle
\stackrel{def}{=}\lim_{\alpha \to n+2s}
\langle \kappa_{\alpha}(x),\phi \rangle
=\lim_{\alpha \to n+2s}
\int_{\bR^n}\frac{|x|^{\alpha-n}}{\gamma_n(\alpha)}\phi(x)\,d^nx
$$
$$
=\lim_{\beta \to 0}
\int_{\bR^n}\frac{|x|^{2s+\beta}-|x|^{2r}}{\gamma_n(n+2s+\beta)}
\phi(x)\,d^nx
=-\lim_{\beta \to 0}
\int_{\bR^n}|x|^{2s}\frac{|x|^{\beta}-1}{\beta \gamma_n(n+2s)}\phi(x)\,d^nx
$$
\begin{equation}
\label{17.1}
\qquad\qquad
=-\int_{\bR^n}\frac{|x|^{2s}\log|x|}{\gamma_n(n+2s)}\phi(x)\,d^nx,
\quad \forall \, \phi\in \Phi(\bR^n),
\end{equation}
where $|\alpha-n-2s|<2$,
\begin{equation}
\label{17.2}
\gamma_n(n+2s)=(-1)^{s}2^{n+2s-1}\pi^{\frac{n}{2}}
s!\Gamma\Big(\frac{n}{2}+s\Big), \quad s=0,1,2,\dots.
\end{equation}
Thus,
\begin{equation}
\label{18}
\kappa_{n+2s}(x)\stackrel{def}{=}\lim_{\alpha \to n+2s}\kappa_{\alpha}(x)
=-\frac{|x|^{2s}\log|x|}{\gamma_n(n+2s)}, \quad s=0,1,2,\dots.
\end{equation}
Formulas (\ref{17}), (\ref{18}) for the {\it Riesz kernel\/}
were constructed in~\cite{Liz3} (see
also~\cite[Lemma~2.13.]{Sam3},~\cite[Lemma~25.2.]{Sam-Kil-Mar}).

If $n=1$ we have
\begin{equation}
\label{19}
\kappa_{\alpha}(x)=\left\{
\begin{array}{lcl}
\displaystyle
\frac{|x|^{\alpha-1}\Gamma(\frac{1-\alpha}{2})}
{2^{\alpha}\sqrt{\pi}\Gamma(\frac{\alpha}{2})}
=\frac{|x|^{\alpha-1}}{2\Gamma(\alpha)\cos(\frac{\pi\alpha}{2})},
&& \alpha \ne -2s, \, \alpha \ne 1+2s, \\
\displaystyle
(-1)^{s+1}\frac{|x|^{2s}\log|x|}{\pi(2s)!},
&& \alpha=1+2s, \\
\displaystyle
(-1)^s\delta^{(2s)}(x),
&& \alpha=-2r, \, s\in \bN_0, \\
\end{array}
\right.
\end{equation}

Easy calculations show that if $\alpha \ne n+2r$ then
the Riesz kernel $\kappa_{\alpha}(x)$ is a {\it homogeneous\/}
distribution of degree~$\alpha-n$, and if $\alpha=n+2s$ then
the Riesz kernel is an {\it associated homogeneous\/} distribution
of degree $2s$ and order $1$, \ $s\in {\bN}_0$ (see
Definitions~\ref{de2},~\ref{de2.1}).

According to~\cite[Ch.II,\S3.3.,(2)]{Gel-S},~\cite[Lemma~2.13.]{Sam3},
~\cite[Lemma~25.2.]{Sam-Kil-Mar},
\begin{equation}
\label{20}
F[\kappa_{\alpha}(x)](\xi)=|x|^{-\alpha},
\end{equation}
where for $\alpha=n+2s$ the right-hand side of the last relation is understood
as the principal value given by (\ref{12}).

With the help of (\ref{6}), (\ref{7}), and (\ref{20}), we obtain
$$
\kappa_{\alpha}(x)*\kappa_{\beta}(x)=\kappa_{\alpha+\beta}(x),
\quad Re\,\alpha, \ Re\,\beta <0,
$$
where $\alpha, \, \beta, \, \alpha+\beta\ne -2s$, \
$s=0,1,2,\dots$.
Next, by analytic continuation of the left-hand and right-hand sides
of the last formula with respect to $\alpha$, and taking into account
formula (\ref{17}), we define the relation
\begin{equation}
\label{21}
\kappa_{\alpha}(x)*\kappa_{\beta}(x)=\kappa_{\alpha+\beta}(x),
\end{equation}
in the sense of distribution from ${\cS}'(\bR^n)$, where
$\alpha, \, \beta, \, \alpha+\beta \ne n+2s$, \ $s=0,1,2,\dots$.
Taking into account formula (\ref{18}), it is easy to see that
\begin{equation}
\label{21.1}
\kappa_{\alpha}(x)*\kappa_{\beta}(x)=\kappa_{\alpha+\beta}(x),
\quad \alpha, \beta \in \bC,
\end{equation}
in the sense of distribution from $\Phi'(\bR^n)$.

Let $\alpha=(\alpha_1,\dots,\alpha_n)$, $\alpha_j\in \bC$, $j=1,2,\dots$.
We denote by
\begin{equation}
\label{21.3}
f_{\alpha}(x)=\kappa_{\alpha_1}(x_1)\times\cdots\times
\kappa_{\alpha_n}(x_n),
\end{equation}
the {\it multi-Riesz kernel\/}, where the one-dimensional
Riesz kernel $\kappa_{\alpha_j}(x_j)$, $j=1,\dots,n$ is
defined by (\ref{19}).

If $\alpha_j \ne 1+2r_j$, $r_j\in \bN_0$, \ $j=1,2,\dots$
then the Riesz kernel
$$
f_{\alpha}(x)=\frac{\times_{j=1}^{n}|x_j|^{\alpha_j-1}}
{2^{n}\prod_{j=1}^{n}\Gamma(\alpha_j)\cos(\frac{\pi\alpha_j}{2})}
$$
is a {\it homogeneous\/} distribution of degree~$|\alpha|-n$
(see Definition~\ref{de2}).

If $\alpha_j=1+2r_j$, $r_j\in \bN_0$, \ $j=1,2,\dots,k$, \
$\alpha_s \ne 1+2r_s$, $r_s\in \bN_0$, \ $s=k+1,\dots,n$ then
$$
f_{\alpha}(x)
=\frac{(-1)^{r_1+\cdots+r_k+k}\times_{j=1}^{k}|x_j|^{2r_j}\log|x_j|}
{\pi^{k}(2r_1)!\cdots(2r_k)!}
\qquad\qquad\qquad\qquad\qquad
$$
\begin{equation}
\label{21.4}
\qquad\qquad
\times\frac{\times_{j=k+1}^{n}|x_j|^{\alpha_j-1}}
{2^{n-k}\Gamma(\alpha_{k+1})\cos(\frac{\pi\alpha_{k+1}}{2})\cdots
\Gamma(\alpha_{n})\cos(\frac{\pi\alpha_{n}}{2})}.
\end{equation}
Thus if among all $\alpha_1,\dots,\alpha_n$ there are $k$ pieces such
that $=1+2r$ and $n-k$ pieces such that $\ne 1+2r$, \ $r\in \bN_0$,
then the Riesz kernel $f_{\alpha}(x)$ is an {\it associated homogeneous\/}
distribution of degree~$|\alpha|-n$ and order $k$, \ $k=1,\dots,n$
(see Definition~\ref{de2.1}).

Taking into account the above calculation, it is easy to see that
\begin{equation}
\label{21.5}
f_{\alpha}(x)*f_{\beta}(x)=f_{\alpha+\beta}(x),
\quad \alpha, \beta \in \bC^n,
\end{equation}
in the sense of distribution from $\Phi'_{\times}(\bR^n)$.

\subsection{The Riesz fractional operator.}\label{s2.3}
Define the Riesz fractional operator
$D^{\alpha}:\phi \mapsto D^{\alpha}\phi$ on the Lizorkin space
$\Phi(\bR^n)$ as a convolution
$$
\big(D^{\alpha}\phi\big)(x)
\stackrel{def}{=}(-\Delta)^{\alpha/2}\phi(x)
\qquad\qquad\qquad\qquad\qquad\qquad\qquad\qquad\qquad
$$
\begin{equation}
\label{22}
\quad
=\kappa_{-\alpha}(x)*\phi(x)
=\Bigl\langle \frac{|\xi|^{-\alpha-n}}{\gamma_n(\alpha)},
\phi(x-\xi)\Bigr\rangle,
\quad x\in \bR^n,
\end{equation}
where $\phi \in \Phi(\bR^n)$.

It is known that in the general case, for $\varphi \in {\cS}(\bR^n)$,
the function $(D^{\alpha}\varphi)(x) \not\in {\cS}(\bR^n)$.
However, the following assertion holds.

\begin{Lemma}
\label{lem2}
{\rm (~\cite[Theorem~2.16.]{Sam3},
~\cite[Theorem~25.1.]{Sam-Kil-Mar})}
The Lizorkin space of test functions $\Phi(\bR^n)$ is invariant
under the Riesz fractional operator $D^{\alpha}$ and
$D^{\alpha}(\Phi(\bR^n))=\Phi(\bR^n)$.
\end{Lemma}

\begin{proof}
Indeed, according to (\ref{22}), (\ref{20}), (\ref{7}),
$$
F[D^{\alpha}\phi](\xi)=|\xi|^{-\alpha}F[\phi](\xi),
\quad \phi \in \Phi(\bR^n).
$$
Since $F[\phi](\xi)\in \Psi(\bR^n)$ and
$|\xi|^{-\alpha}F[\phi](\xi)\in \Psi(\bR^n)$ then
$D^{\alpha}\phi \in \Phi(\bR^n)$, i.e.,
$D^{\alpha}(\Phi(\bR^n))\subset \Phi(\bR^n)$. Moreover,
any function from $\Psi(\bR^n)$ can be represented as
$\psi(\xi)=|\xi|^{\alpha}\psi_1(\xi)$, $\psi_1 \in \Psi(\bR^n)$.
This implies that $D^{\alpha}(\Phi(\bR^n))=\Phi(\bR^n)$.
\end{proof}

It is clear that~\cite[(25.2)]{Sam-Kil-Mar}
\begin{equation}
\label{23}
\big(D^{\alpha}\phi\big)(x)=F^{-1}[|\xi|^{\alpha} F[\phi](\xi) ](x),
\quad \phi \in \Phi(\bR^n).
\end{equation}

The operator $D^{\alpha}$ is called the operator of (fractional)
partial differentiation of order $\alpha$, for
$\alpha>0$, and the operator of (fractional) partial integration
of order $\alpha$, for $\alpha <0$; \ $D^{0}$ is the identity operator.

In particular,
$$
\begin{array}{rcl}
\displaystyle
D^{-n-2r}\phi&=&(-\Delta)^{-n/2-r}\phi
\displaystyle
=\frac{(-1)^{r+1}|x|^{2r}\log|x|}
{2^{n+2r-1}\pi^{n/2}r!\Gamma(\frac{n}{2}+r)}*\phi, \\
\displaystyle
D^{2r}\phi&=&(-\Delta)^{r}\phi, \quad \phi \in \Phi(\bR^n),
\qquad r=0,1,2,\dots. \\
\end{array}
$$
Note that definition (\ref{22})
$D^{\alpha}\phi\stackrel{def}{=}(-\Delta)^{\alpha/2}\phi$ is
introduced in view of the last relation.

According to formulas (\ref{6}), (\ref{22}), we define the Riesz
fractional operator $D^{\alpha}f$, $\alpha \in \bC$ of a distribution
$f\in \Phi'(\bR^n)$ by the relation
\begin{equation}
\label{24}
\langle D^{\alpha}f,\phi\rangle\stackrel{def}{=}
\langle f(-\Delta)^{\alpha/2}f,\phi\rangle
=\langle f, D^{\alpha}\phi\rangle,
\quad \forall \, \phi\in  \Phi(\bR^n).
\end{equation}

It is clear that $D^{\alpha}(\Phi'(\bR^n))=\Phi'(\bR^n)$. Moreover,
the family of operators $D^{\alpha}$, $\alpha\in \bC$ forms an
Abelian group on the space $\Phi'(\bR^n)$: if $f \in \Phi'(\bR^n)$ then
\begin{equation}
\label{24.1}
\begin{array}{rcl}
\displaystyle
D^{\alpha}D^{\beta}f&=&D^{\beta}D^{\alpha}f=D^{\alpha+\beta}f, \\
\displaystyle
D^{\alpha}D^{-\alpha}f&=&f,
\qquad \alpha, \, \beta \in \bC. \\
\end{array}
\end{equation}

\subsection{The multi-Riesz fractional operator.}\label{s2.3-1}
Define the multi-Riesz fractional operator
$D^{\alpha}_{\times}: \phi(x) \to D^{\alpha}_{\times}\phi(x)$
on the Lizorkin space $\Phi_{\times}(\bR^n)$ as the convolution
\begin{equation}
\label{24.2}
\Big(D^{\alpha}_{\times}\phi\Big)(x)\stackrel{def}{=}f_{-\alpha}(x)*\phi(x),
\quad \phi\in \Phi_{\times}(\bR^n),
\end{equation}
where the {\it multi-Riesz kernel\/} $f_{-\alpha}(x)$ is given by formula
(\ref{21.3}). Here
$D^{\alpha}_{\times}=D^{\alpha_1}_{x_1}\cdots D^{\alpha_n}_{x_n}$, where
$D^{\alpha_j}_{x_j}=f_{-\alpha_j}(x_j)*$, $j=1,2,\dots,n$.

It is easy to verify that for the operator $D^{\alpha}_{\times}$
an analog of Lemma~\ref{lem2} holds.
\begin{Lemma}
\label{lem2.1}
The Lizorkin space of test functions $\Phi_{\times}(\bR^n)$ is invariant
under the Riesz fractional operator $D^{\alpha}_{\times}$ and
$D^{\alpha}_{\times}(\Phi_{\times}(\bR^n))=\Phi_{\times}(\bR^n)$.
\end{Lemma}

Analogously to (\ref{24}), if $f\in \Phi'_{\times}(\bR^n)$  then
\begin{equation}
\label{24*}
\langle D^{\alpha}_{\times}f,\phi\rangle\stackrel{def}{=}
=\langle f, D^{\alpha}_{\times}\phi\rangle,
\quad \forall \, \phi\in  \Phi_{\times}(\bR^n),
\quad \alpha \in \bC^n.
\end{equation}

The family of operators $D^{\alpha}_{\times}$, $\alpha \in \bC^n$
forms an Abelian group on the space $\Phi'_{\times}(\bR^n)$.

\subsection{Quasi-asymptotics.}\label{s2.4}
Recall the definitions of a {\it quasi-asymptotics\/} \cite{D-Zav1},
\cite{Vl-D-Zav}.

\begin{Definition}
\label{de4.1} \rm
(~\cite[\S3.2.]{Vl-D-Zav}) A positive continuous real-valued function
$\rho(a)$, $a\in \bR$ such that for any $a>0$ there exists the
following limit
$$
\lim_{t \to \infty}\frac{\rho(ta)}{\rho(t)}=C(a)
$$
is called an {\it automodel {\rm(}or regular varying{\rm)}\/} function.
\end{Definition}

It is easy to see that the function $C(a)$ satisfies the functional
equation $C(ab)=C(a)C(b)$, \ $a,b >0$. It is well known that the
solution of this equation is the following:
\begin{equation}
\label{25}
C(a)=a^{\alpha}, \qquad \alpha\in \bR.
\end{equation}
In this case we say that an {\it automodel\/} function $\rho(a)$
has the degree $\alpha$.

For example, the functions $t^{\alpha}$, \ $t^{\alpha}\log^{m}t$,
$m\in \bN$ ($t>0$) are {\it automodel\/} of degree $\alpha$.

\begin{Definition}
\label{de4} \rm
Let $f\in {\cS}'(\bR^n)$. If there exists an
{\it automodel\/} function $\rho(t)$, $t>0$ of degree $\alpha$ such that
$$
\frac{f(tx)}{\rho(t)} \to g(x)\not\equiv 0, \quad t \to \infty,
\quad \text{in} \quad {\cS}'(\bR^n).
$$
then we say that the distribution $f$ has the {\it quasi-asymptotics\/}
$g(x)$ of degree $\alpha$ at infinity with respect to $\rho(t)$,
and write
$$
f(x) \stackrel{{\cS}'}{\sim} g(x), \quad |x| \to \infty \ \big(\rho(t)\big).
$$
If for any $\alpha$ we have
$$
\frac{f(tx)}{t^{\alpha}} \to 0, \quad t \to \infty,
\quad \text{in} \quad {\cS}'(\bR^n)
$$
then we say that the distribution $f$ has a {\it quasi-asymptotics\/}
of degree $-\infty$ at infinity and write $f(x) \stackrel{{\cS}'}{\sim} 0$, \
$|x| \to \infty$.
\end{Definition}

\begin{Lemma}
\label{lem1}
{\rm (~\cite{D-Zav1},~\cite[\S3.2.]{Vl-D-Zav})}
Let $f\in {\cS}'(\bR^n)$. If $f(x) \stackrel{{\cS}'}{\sim} g(x)\not\equiv 0$,
as $|x| \to \infty$ with respect to the {\it automodel\/} function
$\rho(t)$ of degree $\alpha$ then $g(x)$ is a homogeneous distribution
of degree $\alpha$.
\end{Lemma}

If $n=1$, the results from~\cite[Ch.I,\S3.11.]{Gel-S} and
Lemma~\ref{lem1} imply that
\begin{equation}
\label{18.2}
g(x)=\left\{
\begin{array}{lcr}
C_{1}x_{+}^{\alpha}+C_{2}x_{-}^{\alpha}, &&\quad \alpha \ne -k, \\
C_{1}P(x^{-k})+C_{2}\delta^{(k-1)}(x),      &&\quad \alpha=-k, \\
\end{array}
\right.
\end{equation}
where $C_1,C_2$ are constants.

Here the distributions $x_{\pm}^{\alpha}$, $\alpha \ne -k$, \ $k\in \bN$
are defined by the following relations~\cite[Ch.I,\S3.2]{Gel-S}:
if $Re \alpha>-m-1$, $\alpha\ne -1,-2,\dots,-m$, \ $m\in \bN_0$ then
$$
\bigl\langle x_{+}^{\alpha},\varphi(x) \bigr\rangle
\stackrel{def}{=}\int_{0}^{1}x^{\alpha}\bigg(\varphi(x)
-\sum_{j=0}^{m-1}\frac{x^{j}}{j!}\varphi^{(j)}(0)\bigg)\,dx
\qquad\qquad\qquad\qquad
$$
\begin{equation}
\label{18.3}
\qquad\qquad\qquad
+\int_{1}^{\infty}x^{\alpha}\varphi(x)\,dx
+\sum_{j=0}^{m-1}\frac{\varphi^{(j)}(0)}{j!(\lambda+j+1)},
\end{equation}
\begin{equation}
\label{18.4}
\bigl\langle x_{-}^{\alpha},\varphi(x)\bigr\rangle
\stackrel{def}{=}\bigl\langle x_{+}^{\alpha},\varphi(-x) \bigr\rangle,
\end{equation}
for all $\varphi\in {\cS}(\bR)$.
Using (\ref{18.3}), (\ref{18.4}), one can introduce the
distributions~\cite[Ch.I,\S3.3]{Gel-S}:
\begin{equation}
\label{18.5}
\begin{array}{rcl}
\displaystyle
|x|^{\alpha}&\stackrel{def}{=}&x_{+}^{\alpha}+x_{-}^{\alpha},
\quad \alpha \ne -2k+1, \\
\displaystyle
|x|^{\alpha}{\rm sign}x&\stackrel{def}{=}&x_{+}^{\alpha}-x_{-}^{\alpha},
\quad \alpha \ne -2k, \\
\end{array}
\end{equation}
where $k\in \bN$. For the other $k$ these distributions are well-defined.
The principal value of the functions~$x^{-2k}$ and $x^{-2k+1}$
are defined as
\begin{equation}
\label{18.6}
\begin{array}{rcl}
\displaystyle
P(x^{-2k+1})&\stackrel{def}{=}&|x|^{-2k+1}{\rm sign\,}x, \\
\displaystyle
P(x^{-2k})&\stackrel{def}{=}&|x|^{-2k}, \quad k\in \bN, \\
\end{array}
\end{equation}
respectively.

\begin{Definition}
\label{de5} \rm
Let $f\in {\cS}'(\bR^n)$. If there exists an {\it automodel\/} function
$\rho(t)$, $t>0$ of degree $\alpha$ such that
$$
\frac{f(\frac{x}{t})}{\rho(t)} \to g(x)\not\equiv 0, \quad t \to \infty,
\quad \text{in} \quad {\cS}'(\bR^n)
$$
then we say that the distribution $f$ has a {\it quasi-asymptotics\/}
$g(x)$ of degree $-\alpha$ at zero with respect to $\rho(t)$, and write
$$
f(x) \stackrel{{\cS}'}{\sim} g(x), \quad |x| \to 0 \ \big(\rho(t)\big).
$$
If for any $\alpha$ we have
$$
\frac{f(\frac{x}{t})}{t^{\alpha}} \to 0, \quad t \to \infty,
\quad \text{in} \quad {\cS}'(\bR^n)
$$
then we say that the distribution $f$ has a {\it quasi-asymptotics\/}
of degree $-\infty$ at zero, and write $f(x) \stackrel{{\cS}'}{\sim} 0$, \
$|x| \to 0$.
\end{Definition}

For the case of distributions from $\Phi'(\bR^n)$ and
$\Phi_{\times}'(\bR^n)$ Definitions~\ref{de4},~\ref{de5}
and Lemma~\ref{lem1} are formulated word for word.

\section{The Tauberian type theorems}
\label{s3}

\begin{Theorem}
\label{th5}
A distribution $f\in {\cS}'(\bR^n)$ has a quasi-asymptotics of degree
$\alpha$ at infinity with respect to the automodel function $\rho(t)$,
$t>0$, if and only if its Fourier transform has a quasi-asymptotics of
degree $-\alpha-n$ at zero with respect to the automodel function $t^n\rho(t)$.
\end{Theorem}

\begin{proof}
Since $F[f(x)](\xi/t)=t^nF[f(tx)](\xi)$, \ $x,\xi\in \bR^n$, \ $t>0$,
we have
$$
\bigl\langle F[f(x)](\xi/t),\varphi(\xi)\bigr\rangle
=t^n\bigl\langle F[f(tx)](\xi),\varphi(\xi)\bigr\rangle
=t^n\bigl\langle f(tx),F[\varphi(\xi)](x)\bigr\rangle,
$$
$\varphi \in {\cS}(\bR^n)$.
Thus
$$
\lim_{t \to \infty}
\Bigl\langle \frac{F[f(x)](\xi/t)}{t^n\rho(t)},\varphi(\xi)\Bigr\rangle
=\lim_{t \to \infty}
\Bigl\langle \frac{f(tx)}{\rho(t)},F[\varphi(\xi)](x)\Bigr\rangle,
\quad \forall \, \varphi \in {\cS}(\bR^n).
$$
The last relation implies that
$f(x) \stackrel{{\cS}'}{\sim} g(x)$, \ $|x| \to \infty$ \
$\big(\rho(t)\big)$, i.e.,
$$
\lim_{t \to \infty}
\Bigl\langle \frac{f(tx)}{\rho(t)}, \varphi(x)\Bigr\rangle
=\langle g(x), \varphi(x)\rangle, \quad \forall \, \varphi \in {\cS}(\bR^n),
$$
if and only if
$F[f(x)](\xi) \stackrel{{\cS}'}{\sim} F[g(x)](\xi)$, \ $|\xi| \to 0$ \
$\big(t^n\rho(t)\big)$, i.e.,
$$
\lim_{t \to \infty}
\Bigl\langle \frac{F[f(x)](\xi/t)}{t^n\rho(t)},\varphi(\xi)\Bigr\rangle
=\bigl\langle F[g(x)](\xi),\varphi(\xi)\bigr\rangle.
$$
\end{proof}

\begin{Theorem}
\label{th6}
A distribution $f\in {\cS}'(\bR)$ has a quasi-asymptotics of
degree $\alpha$ at infinity, i.e.,
$$
f(x) \stackrel{{\cS}'}{\sim}
g(x)=\left\{
\begin{array}{lcr}
C_{1}x_{+}^{\alpha}+C_{2}x_{-}^{\alpha}, &&\quad \alpha \ne -k, \\
C_{1}P(x^{-k})+C_{2}\delta^{(k-1)}(x),   &&\quad \alpha=-k, \\
\end{array}
\right.
\quad |x| \to \infty,
$$
if and only if its Fourier transform $F[f]$ has a
quasi-asymptotics of degree $-\alpha-1$ at zero, i.e.,
$$
F[f(x)](\xi)\stackrel{{\cS}'}{\sim}
F[g(x)](\xi)=\left\{
\begin{array}{lcr}
\Gamma(\alpha+1)\big(B_1\xi_{+}^{-\alpha-1}+B_2\xi_{-}^{-\alpha-1}\big), \\
C_{1}\frac{\pi i^{k}}{(k-1)!}\xi^{k-1}{\rm sign}\xi+C_{2}(-i\xi)^{k-1}, \\
\end{array}
\right.
\quad |\xi| \to 0,
$$
where $C_1$, $C_2$ are constants, and
$B_1=C_1e^{i(\alpha+1)\pi/2}+C_2e^{-i(\alpha+1)\pi/2}$, \
$B_2=C_1e^{-i(\alpha+1)\pi/2}+C_2e^{i(\alpha+1)\pi/2}$, \
$k\in \bN$.
\end{Theorem}

The proof of Theorem~\ref{th6} follows from
Theorem~\ref{th5}, formula (\ref{18.2}), and formulas
from~\cite[Ch.II,\S2.3.,Ch.I,\S3.6.]{Gel-S}.

\begin{Theorem}
\label{th7}
Let $f \in \Phi'(\bR^n)$. Then
$$
f(x) \stackrel{\Phi'}{\sim} g(x), \quad |x| \to \infty
\quad \big(\rho(t)\big)
$$
if and only if
$$
D^{\beta}f(x) \stackrel{\Phi'}{\sim} D^{\beta}g(x), \quad
|x| \to \infty \quad \big(t^{-\beta}\rho(t)\big),
$$
where $\beta\in \bC$.
\end{Theorem}

\begin{proof}
Let $\beta \ne -n-2r$, $r\in \bN_0$. Since the {\it Riesz kernel\/}
(\ref{15}), (\ref{17}) is a {\it homogeneous\/} distribution of
degree~$\alpha-n$, according to Lemma~\ref{lem2} and formulas
(\ref{24}), (\ref{22}), (\ref{6}), we have
$$
\bigl\langle \big(D^{\beta}f\big)(tx),\phi(x)\bigr\rangle
=\bigl\langle \big(f*\kappa_{-\beta}\big)(tx), \phi(x)\bigr\rangle
\qquad\qquad\qquad\qquad\qquad\qquad
$$
$$
\qquad\qquad
=t^{-n}\Bigl\langle f(x),
\Bigl\langle \kappa_{-\beta}(y), \phi\Big(\frac{x+y}{t}\Big)
\Bigr\rangle\Bigr\rangle
=t^{n}
\bigl\langle f(tx),
\bigl\langle \kappa_{-\beta}(ty), \phi(x+y)\bigr\rangle \bigr\rangle
$$
\begin{equation}
\label{41}
=t^{-\beta}
\bigl\langle f(tx),
\bigl\langle \kappa_{-\beta}(y), \phi(x+y)\bigr\rangle \bigr\rangle
=t^{-\beta}
\bigl\langle f(tx), \big(D^{\beta}\phi\big)(x)\bigr\rangle,
\end{equation}
for all $\phi \in \Phi(\bR^n)$. Thus
$$
\Bigl\langle \frac{\big(D^{\beta}f\big)(tx)}{t^{-\beta}\rho(t)},
\phi(x)\Bigr\rangle
=\Bigl\langle \frac{f(tx)}{\rho(t)},\big(D^{\beta}\phi\big)(x)\Bigr\rangle.
$$

Taking into account that the Lizorkin space of test functions
$\Phi(\bR^n)$ is invariant under the Riesz fractional operator
$D^{\beta}$ and passing to the limit in the above relation, as
$t \to \infty$, we obtain
$$
\lim_{t \to \infty}
\Bigl\langle \frac{\big(D^{\beta}f\big)(tx)}{t^{-\beta}\rho(t)},
\phi(x)\Bigr\rangle
=\lim_{t \to \infty}
\Bigl\langle \frac{f(tx)}{\rho(t)},\big(D^{\beta}\phi\big)(x)\Bigr\rangle
=\bigl\langle D^{\beta}g(x),\phi(x)\bigr\rangle.
$$
That is, $\lim_{t \to \infty}
\frac{\big(D^{\beta}f\big)(tx)}{t^{-\beta}\rho(t)}=D^{\beta}g(x)$
in $\Phi'(\bR^n)$ if and only if
$\lim_{t \to \infty}\frac{f(tx)}{\rho(t)}=g(x)$ in
$\Phi'(\bR^n)$.
Thus this case of the theorem is proved.

Next, consider the case $\beta=-n-2r$, $r\in {\bN}_0$. Since
the {\it Riesz kernel\/} (\ref{18}) is an {\it associated homogeneous\/}
distribution of degree $2r$ and order $1$, consequently,
$$
\kappa_{n+2r}(ty)=-t^{2r}\frac{|y|^{2r}\log|y|}{\gamma_n(n+2r)}
-t^{2r}\log t\frac{|y|^{2r}}{\gamma_n(n+2r)}, \quad r=0,1,2,\dots.
$$
for all $t>0$. In view of (\ref{8}),
$\langle |y|^{2r}, \phi(x+y)\bigr\rangle=0$, and we have
$$
\bigl\langle \kappa_{n+2r}(ty), \phi(x+y)\bigr\rangle
=t^{2r}\bigl\langle \kappa_{n+2r}(y), \phi(x+y)\bigr\rangle
=t^{2r}\big(D^{-n-2r}\phi\big)(x).
$$
Repeating the above calculations almost word for word, we prove this case
of the theorem.
\end{proof}

\begin{Theorem}
\label{th8}
Let $f \in \Phi_{\times}'(\bR^n)$. Then
$$
f(x) \stackrel{\Phi_{\times}'}{\sim} g(x), \quad |x| \to \infty
\quad \big(\rho(t)\big)
$$
if and only if
$$
D^{\beta}_{\times}f(x) \stackrel{\Phi_{\times}'}{\sim} D^{\beta}_{\times}g(x),
\quad |x| \to \infty \quad \big(t^{|-\beta|}\rho(t)\big),
$$
where $\beta=(\beta_1,\dots,\beta_n)\in \bC^n$, \
$|\beta|=\beta_1+\cdots+\beta_n$.
\end{Theorem}

\begin{proof}
Let $\beta_j \ne -1-2r_j$, $r_j\in \bN_0$, \ $j=1,\dots,n$.
In this case the Riesz kernel $f_{-\beta}(x)$ is a
{\it homogeneous\/} distribution of degree~$|-\beta|-n$.
Using Lemma~\ref{lem2.1} and formulas (\ref{24*}), (\ref{24.2}),
(\ref{6}), we obtain
$$
\bigl\langle \big(D^{\beta}_{\times}f\big)(tx),\phi(x)\bigr\rangle
=\bigl\langle \big(f*f_{-\beta}\big)(tx),\phi(x)\bigr\rangle
\qquad\qquad\qquad\qquad\qquad\qquad
$$
$$
\qquad
=t^{-n}
\Bigl\langle f(x),\Bigl\langle f_{-\beta}(y),\phi\Big(\frac{x+y}{t}\Big)
\Bigr\rangle\Bigr\rangle
=t^{n}
\bigl\langle f(tx),
\bigl\langle f_{-\beta}(ty),\phi(x+y)\bigr\rangle \bigr\rangle
$$
$$
=t^{|-\beta|}
\bigl\langle f(tx),
\bigl\langle f_{-\beta}(y),\phi(x+y)\bigr\rangle \bigr\rangle,
=t^{|-\beta|}
\bigl\langle f(tx),\big(D^{\beta}_{\times}\phi\big)(x)\bigr\rangle,
\quad
$$
for all $\phi \in \Phi_{\times}(\bR^n)$.
Thus
$$
\Bigl\langle \frac{\big(D^{\beta}_{\times}f\big)(tx)}{t^{|-\beta|}\rho(t)},
\phi(x)\Bigr\rangle
=\Bigl\langle \frac{f(tx)}{\rho(t)},
\big(D^{\beta}_{\times}\phi\big)(x)\Bigr\rangle,
$$
and, consequently, $\lim_{t\to \infty}
\frac{(D^{\beta}_{\times}f)(tx)}{t^{|-\beta|}\rho(t)}=D^{\beta}_{\times}g(x)$
if and only if $\lim_{t\to \infty}\frac{f(tx)}{\rho(t)}=g(x)$ in
$\Phi_{\times}'(\bR^n)$.
Thus this case of the theorem is proved.

Consider the case when among all $\beta_1,\dots,\beta_n$ there are $k$
pieces such that $=-1-2r$ and $n-k$ pieces such that $\ne -1-2r$, \
$r\in \bN_0$. In this case, according to Definition~\ref{de2.1},
the Riesz kernel $f_{-\beta}(x)$ is an {\it associated homogeneous\/}
distribution of degree~$|-\beta|-n$ and order $k$, \ $k=1,\dots,n$.

Let $\beta_j=-1-2r_j$, $r_j\in \bN_0$, \ $j=1,\dots,k$; \
$\beta_s \ne -1-2r_s$, $r_s\in \bN_0$, \ $s=k+1,\dots,n$.
Denote
$$
A=(-1)^{r_1+\cdots+r_k+k}\pi^{k}(2r_1)!\cdots(2r_k)!
\qquad\qquad\qquad\qquad\qquad\qquad
$$
$$
\qquad\qquad\qquad\qquad
\times
2^{n-k}\Gamma(-\beta_{k+1})\cos(\frac{\pi\beta_{k+1}}{2})\cdots
\Gamma(-\beta_{n})\cos(\frac{\pi\beta_{n}}{2}).
$$
Then according to (\ref{21.4}),
$$
f_{-\beta}(ty)=\frac{1}{A}t^{|-\beta|-n}
|y_{1}|^{2r_1}\times\cdots\times |y_k|^{2r_k}
\times|y_{k+1}|^{-\beta_{k+1}-1}\times\cdots\times |y_n|^{-\beta_{n}-1}
$$
$$
\qquad\qquad\qquad\qquad\qquad\qquad
\times(\log|y_1|+\log|t|)\times\cdots\times(\log|y_k|+\log|t|)
$$
$$
=t^{|-\beta|-n}f_{-\beta}(y)
\qquad\qquad\qquad\qquad\qquad\qquad\qquad\qquad\qquad
$$
$$
+\frac{1}{A}|y_{1}|^{2r_1}\times\cdots\times |y_k|^{2r_k}
\times|y_{k+1}|^{-\beta_{k+1}-1}\times\cdots\times |y_n|^{-\beta_{n}-1}
$$
$$
\times
\bigg(\Big(\log|y_2|\times\cdots\times\log|y_k|
+\cdots+\log|y_1|\times\cdots\times\log|y_{k-1}|\Big)\log t
$$
\begin{equation}
\label{42}
\qquad
+\cdots+\Big(\log|y_1|+\cdots+\log|y_k|\Big)\log^{k-1}|t|
+\log^{k}t\bigg).
\end{equation}

It is easy to verify that in view of characterization (\ref{8*}),
\begin{equation}
\label{42*}
\bigl\langle f_{-\beta}(ty), \phi(x+y)\bigr\rangle
=t^{|-\beta|-n}\bigl\langle f_{-\beta}(y), \phi(x+y)\bigr\rangle
=t^{|-\beta|-n}\big(D^{\beta}_{\times}\phi\big)(x),
\end{equation}
$\phi \in \Phi_{\times}(\bR^n)$.
For example, taking into account (\ref{8*}), we obtain
$$
\Bigl\langle
\times_{j=2}^{k}(x_{j}-y_{j})^{2r_{j}}\log|x_j-y_j|
\times_{i=k+1}^{n}|x_{i}-y_{i}|^{-\beta_{i}-1},
\qquad\qquad\qquad\qquad
$$
$$
\qquad\qquad\qquad\qquad
\int_{\bR}(x_{1}-y_{1})^{2r_{1}}\phi(y_1,y_2,\dots,y_n)\,dy_1 \Bigr\rangle=0,
$$
for all $\phi \in \Phi_{\times}(\bR^n)$.
In a similar way, one can prove that all terms in (\ref{42}),
with the exception of $t^{|-\beta|-n}f_{-\beta}(y)$
{\it do not make a contribution\/} to the functional
$\langle f_{-\beta}(ty),\phi(x+y)\rangle$, where
$\beta_j=-1-2r_j$, $r_j\in \bN_0$, \ $j=1,\dots,k$; \
$\beta_s \ne -1-2r_s$, $r_s\in \bN_0$, \ $s=k+1,\dots,n$.

Thus repeating the above calculations almost word for word and using
(\ref{42*}), we prove this case of the theorem.
\end{proof}

\begin{Theorem}
\label{th9}
A distribution $f\in \Phi'(\bR)$ has an even quasi-asymptotics at
infinity with respect to an automodel function $\rho(t)$ of degree
$\alpha$ if and only if there exists a positive integer $N>-\alpha$
such that
$$
\lim_{|x| \to \infty}\frac{D^{-N}f(x)}{|x|^{N}\rho(x)}=A \ne 0,
$$
i.e., the fractional primitive $D^{-N}f(x)$ of
order $N$ has an asymptotics of degree $\alpha+N$ at infinity
{\rm(}understood in the usual sense{\rm)}.
\end{Theorem}

\begin{proof}
Let us prove the necessity.
In view of Theorem~\ref{th7}, setting $\beta=-N$ we have
\begin{equation}
\label{43}
\lim_{t \to \infty}
\Bigl\langle \frac{(D^{-N}f)(tx)}{t^{N}\rho(t)},\phi(x)\Bigr\rangle
=\bigl\langle D^{-N}g_{\alpha}(x),\phi(x)\bigr\rangle
=\bigl\langle \kappa_{N}(x)*g_{\alpha}(x),\phi(x)\bigr\rangle,
\end{equation}
for all $\phi \in \Phi(\bR)$. Here $\kappa_{\beta}(x)$ is
given by (\ref{19}) and $g_{\alpha}(x)$ by (\ref{18.2}).
Since $g(x)$ is even, in view of (\ref{18.2})--(\ref{18.6}),
and (\ref{19}), $g_{\alpha}(x)=C\kappa_{\alpha+1}(x)$,
where $g_{-2k-1}(x)=C\delta^{(2k)}(x)$, \ $k\in \bN_0$.

With the help of formulas (\ref{21.1}), (\ref{19}), we calculate that
$$
\kappa_{N}(x)*g_{\alpha}(x)=C\kappa_{N+\alpha+1}=A|x|^{\alpha+N},
$$
where $A$ is a constant. Thus Eq.~(\ref{43}) can be rewritten as
\begin{equation}
\label{45}
\lim_{t \to \infty}
\Bigl\langle \frac{(D^{-N}f)(tx)}{t^{N}\rho(t)},\phi(x)\Bigr\rangle
=C\bigl\langle \kappa_{N+\alpha+1},\phi(x)\bigr\rangle,
=A\bigl\langle |x|^{\alpha+N},\phi(x)\bigr\rangle,
\end{equation}
for all $\phi \in \Phi(\bR)$. Taking into account that
$N+\alpha>0$, we have
\begin{equation}
\label{46}
\lim_{t \to \infty}
\frac{(D^{-N}f)(tx)}{t^{N}\rho(t)}=A|x|^{\alpha+N}.
\end{equation}

By using Definition~\ref{de4.1} and formula (\ref{25}), relation
(\ref{46}) can be rewritten in the following form
$$
A=\lim_{t \to \infty}
\frac{(D^{-N}f)(tx)}{t^{N}\rho(t)|x|^{\alpha+N}}
\qquad\qquad\qquad\qquad\qquad\qquad\qquad\qquad
$$
\begin{equation}
\label{47}
=\lim_{|tx| \to \infty}\frac{(D^{-N}f)(tx)}{|tx|^{N}\rho(tx)}
\lim_{t \to \infty}\frac{\rho(tx)}{|x|^{\alpha}\rho(t)}
=\lim_{|y| \to \infty}\frac{\big(D^{-N}f\big)(y)}{|y|^{N}\rho(y)}.
\end{equation}

Now we prove the necessity. Relation (\ref{47}) implies (\ref{46}).
The last relation can be rewritten in the weak sense as (\ref{45}).
Next, we rewrite (\ref{45}) in the form (\ref{43}) and using
Theorem~\ref{th7}, prove our assertion.
\end{proof}

\end{document}